\documentclass[english]{article}
\usepackage[T1]{fontenc}
\usepackage[latin9]{inputenc}
\usepackage{geometry}
\usepackage{units}
\usepackage{amsmath}
\usepackage{amssymb}
\usepackage{esint}

\makeatletter
\@ifundefined{definecolor}{\usepackage{color}}{}
\@ifundefined{definecolor}{\usepackage{color}}{}

\usepackage{babel}\usepackage{nicefrac}\usepackage{sagetex}

\usepackage{babel}

\makeatother

\usepackage{babel}
\begin{document}

\title{Computational Aspects of the Combinatorial Nullstellensatz Method}

\author{Edinah K. Gnang%
\thanks{School of Mathematics, Institute for Advanced Study. Email: gnang@ias.edu.%
}}
\maketitle
\begin{abstract}
We discuss here some computational aspects of the Combinatorial Nullstellensatz
argument. Our main result shows that the order of magnitude of the
symmetry group associated with permutations of the variables in algebraic
constraints, determines the performance of algorithms naturally deduced
from Alon's Combinatorial Nullstellensatz arguments. Finally we present
a primal-dual polynomial constructions for certifying the existence
or the non-existence of solutions to combinatorial problems. 
\end{abstract}

\section{Introduction}

It is well-known that systems of polynomial equations over algebraically
closed fields provide concise encodings for classical NP-hard problems
such as the subgraph isomorphism problem. In \cite{A}, Alon presents
the Combinatorial Nullstellensatz method, a general unified algebraic
framework for establishing the existence of solutions to numerous
problems in combinatorics and combinatorial number theory. The Combinatorial
Nullstellensatz argument has recently been subject to intense scrutiny
in the literature. Laso\'{n} in \cite{La} proposed a generalization
of Alon's original formulation of the Combinatorial Nullstellensatz
by weakening the assumption on the degree of the nonvanishing leading
monomial. Furthermore, several variations of the proof of the Combinatorial
Nullstellensatz can be found in the literature \cite{Mi,H,Ko}.

In the concluding remarks of \cite{A} Alon points out that the proofs
presented in \cite{A} are based on non-constructive algebraic arguments
and thus supply no efficient procedure for solving the corresponding
algorithmic problems. Alon further proceeded to raise the fundamental
problem of determining whether or not it is possible to deduce from
such arguments efficient procedures for solving the corresponding
algorithmic problems. Following up on the fundamental problem raised
by Alon, we remark that it is well-known that combinatorial problems
formulated as systems of polynomial equations can be solved using
standard tools in computational algebra such as Gröbner basis \cite{B,CLO}
and closely related methods as presented in \cite{S}. Unfortunately
a precise analysis of the performance of Gröbner basis approaches
in relation to special instances of combinatorial problems remain
unknown. In subsequent work \cite{LMM,LMO,M,LHMO}, the authors follow
up on the fundamental problem raised in \cite{A} by Alon and propose
the Nullstellensatz Linear Algebra algorithmic framework. Their proposed
algorithms, relies on the experimentally-observed low degrees of Hilbert's
Nullstellensatz certificates polynomial encodings of special families
of combinatorial problem instances. The research program developed
in \cite{LMM,LMO,LHMO} follows up on the natural connection between
Hilbert's Nullstellensatz \cite{K} and complexity theory. This connection
was first pointed out by Lovasz in \cite{Lo}. Margulies further developed
this connection in \cite{M}, and established that the minimum-degree
Nullstellensatz certificate for the non-existence of an independent
set of size greater than the largest independent set in the graph
is equal to the size of the largest independent set in the graph.
Finally, in \cite{LMM,LMO,M,LHMO} the authors suggest that algebraic
formulations enable us to exploit sparsity structure of special families
of combinatorial problems. The authors also suggest that symmetries
with respect to permutations of the variables in the algebraic constraints
could yield performance improvements for algorithms which derive Nullstellensatz
certificates. This last suggestion by the authors that exploiting
symmetries are helpful for deriving Nullstellensatz certificates stems
from the well-known role of symmetries for solving combinatorial problems
\cite{SSK,K}.

Our main result establishes that the order of magnitude of the automorphism
group of the constraints with respect to permutations of the variables,
determines the performance of algorithms deduced from Alon's Combinatorial
Nullstellensatz argument for NP-hard problems. We further show that
the Combinatorial Nullstellensatz method yields a natural framework
for a primal-dual certificates for the existence versus the non-existence
of solutions to graph and subgraph isomorphism instances.

\section{Preliminary.}

The Hadamard product of two given column vectors $\mathbf{a}$, $\mathbf{b}\in\mathbb{C}^{n\times1}$
noted $\mathbf{a}\star\mathbf{b}$, corresponds to a column vector
of the same dimensions whose entries correspond to the product of
the corresponding entries of $\mathbf{a}$, and $\mathbf{b}$; we
write 
\begin{equation}
k-\mbox{th entry of }\mathbf{a}\star\mathbf{b}\mbox{ is }a_{k}b_{k}.
\end{equation}
Let us recall here the familiar notation used for the vector product
of $\mathbf{a}$, $\mathbf{b}\in\mathbb{C}^{n\times1}$ with background
matrix the $n\times n$ matrix $\mathbf{M}$. We write 
\begin{equation}
\left\langle \mathbf{a},\:\mathbf{b}\right\rangle _{\mathbf{M}}:=\sum_{0\le k_{0},k_{1}<n}a_{k_{0}}m_{k_{0},k_{1}}b_{k_{1}},
\end{equation}
in particular it follows that 
\begin{equation}
\left\langle \mathbf{a},\:\mathbf{b}\right\rangle :=\left\langle \mathbf{a},\:\mathbf{b}\right\rangle _{\mathbf{I}}=\sum_{0\le k<n}a_{k}b_{k}\:.
\end{equation}
It shall be convenient to adopt the notation convention 
\begin{equation}
\mathbf{a}^{\star^{\alpha}}:=\left(\left(a_{k}\right)^{\alpha}\right)_{0\le k<n}
\end{equation}
Futhermore let $\omega_{n}$ the denote the primitive $n$-th root
of unity expressed by 
\begin{equation}
\omega_{n}=e^{\frac{2\pi i}{n}}
\end{equation}
and
\begin{equation}
\Omega_{n}:=\left\{ \left(\omega_{n}\right)^{k}\right\} _{0\le k<n}.
\end{equation}
The Discrete Fourier Transform matrix $\mathbf{W}$ whose entries
are specified as follows 
\begin{equation}
\mathbf{W}:=\left(w_{uv}=\left(\omega_{n}\right)^{u\cdot v}\right)_{0\le u,v<n}
\end{equation}
is such that 
\begin{equation}
\mathbf{W}\cdot\mathbf{W}^{\dagger}=n\,\mathbf{I}=\mathbf{W}^{\dagger}\cdot\mathbf{W}.
\end{equation}
We shall often denote the set of column vectors of the DFT matrix
$\mathbf{W}$ by the set $\left\{ \mathbf{w}^{\star^{k}}\right\} _{0\le k<n}$,
where
\begin{equation}
\mathbf{w}\,:=\left(w_{k}=\left(\omega_{n}\right)^{k}\right)_{0\le k<n}.
\end{equation}

\section{Overview of the Combinatorial Nullstellensatz}

The Combinatorial Nullstellensatz first presented in \cite{A} by
Alon , corresponds to the following theorem.\\
 \\
\textbf{Theorem} (Combinatorial Nullstellensatz I\textbf{ Alon 1999}):
Let $\mathbb{F}$ be an arbitrary field, and let $f\in\mathbb{F}\left[\mathbf{x}\right]$
and $\left\{ S_{k}\subset\mathbb{F}\right\} _{0\le k<n}$ denote a
collection of non-empty subsets of $\mathbb{F}$ and define $g_{i}\left(x_{i}\right)=\prod_{s\in S_{i}}\left(x_{i}-s\right)$.
If $f$ vanishes over all common zeros of the vector $\mathbf{g}=\left(g_{i}\right)_{0\le i<n}\in\left(\mathbb{F}\left[\mathbf{x}\right]\right)^{n}$
(that is; if $f\left(\mathbf{s}\right)=0$ forall $\mathbf{s}\in S_{0}\times\cdots\times S_{n-1}$
), then there is a vector $\mathbf{h}=\left(h_{i}\right)_{0\le i<n}\in\left(\mathbb{F}\left[\mathbf{x}\right]\right)^{n}$
satisfying the innequalities $\left\{ \deg\left(h_{i}\right)\le\deg\left(f\right)-\deg\left(g_{i}\right)\right\} _{0\le i<n}$
such that 
\begin{equation}
f\left(\mathbf{x}\right)=\left\langle \mathbf{h}\left(\mathbf{x}\right),\mathbf{g}\left(\mathbf{x}\right)\right\rangle :=\sum_{0\le k<n}h_{k}\left(\mathbf{x}\right)\, g_{k}\left(\mathbf{x}\right).
\end{equation}
Moreover, if the polynomial $f$ and the entries of $\mathbf{g}$
lie in $R\left[\mathbf{x}\right]$ for some subring $R$ of $\mathbb{F}$
there are polynomials $\left\{ h_{i}\right\} _{0\le i<n}\subset R\left[\mathbf{x}\right]$
as above.\\
 \\
Consequently we can prove the following :\\
 \\
\textbf{Theorem} (Combinatorial Nullstellensatz II\textbf{ Alon 1999}):
Let $\mathbb{F}$ be an arbitrary field, and let $f\in\mathbb{F}\left[\mathbf{x}\right]$.
Suppose $\deg\left(f\right)$ is $\sum_{0\le i<n}t_{i}$, where each
$t_{i}$ is a nonnegative integer, and suppose the coefficient of
$\mathbf{x}^{\mathbf{t}}$ in $f$ is nonzero. Then if $\left\{ S_{k}\subset\mathbb{F}\right\} _{0\le k<n}$
denotes a collection of non-empty subset of $\mathbb{F}$ with $\left|S_{i}\right|>t_{i}$,
there are $\mathbf{s}\in S_{0}\times\cdots\times S_{n-1}$ so that
\begin{equation}
f\left(\mathbf{s}\right)\ne0.
\end{equation}
 Subsequently the following generalization of the Combinatorial Nullstellensatz
was proposed by Laso\'{n} in \cite{La}.\\
 \\
\textbf{Theorem (}Generalized Combinatorial Nullstellensatz\textbf{
M. Laso\'{n} 2013):} Let $f$ $\in\mathbb{F}\left[\mathbf{x}\right]$.
If $\mathbf{x}^{\mathbf{t}}$ is a nonvanishing monomial in $f$ and
$\mathbf{t}$ is maximal in $\mbox{Supp}(f)$, then for any subsets
$S_{0},\ldots,S_{n-1}$ of $\mathbb{F}$ satisfying $\left\vert S_{i}\right\vert >t_{i}$,
there are $\mathbf{s}\in S_{0}\times\cdots\times S_{n-1}$ so that
$f\left(\mathbf{s}\right)\neq0$.\\
 \\
We present here for the sake of completeness a proof by contradiction
of Alon's Combinatorial Nullstellensatz argument but we remark that
alternative short and possibly more elegant proofs can be found in
\cite{A,Mi,La,Kou,H}.\\
 \\
\textbf{Proof : }We recall that by the Lagrange interpolation formula,
the coefficients of the polynomial $f\in\mathbb{F}\left[\mathbf{x}\right]$
of degree $\sum_{0\le k<n}\left(-1+\left|S_{k}\right|\right)$ are
determined by evalutions of $f$ on the Cartesian product set $S_{0}\times\cdots\times S_{n-1}$.
More explicitly we have 
\[
f\left(\mathbf{x}\right)\equiv\sum_{\mathbf{r}\in S_{0}\times\cdots\times S_{n-1}}f\left(\mathbf{r}\right)\prod_{\left\{ s_{k}\in S_{k}\backslash\left\{ r_{k}\right\} \right\} _{0\le k<n}}\left(\frac{x_{k}-s_{k}}{r_{k}-s_{k}}\right)\mod\left\{ \prod_{s_{i}\in S_{i}}\left(x_{i}-s_{i}\right)\right\} _{0\le i<n},
\]
where the right hand side of the congruence identity corresponds to
the minimal degree polynomial congruent to $f$. By hypothesis we
assumed that the coefficient of the leading monomial $\prod_{0\le k<n}\left(x_{k}\right)^{-1+\left|S_{k}\right|}$
term is non-zero, and incidentally if we have that 
\[
\forall\:\mathbf{r}\in S_{0}\times\cdots\times S_{n-1},\; f\left(\mathbf{r}\right)=0
\]
it would follow that the polynomial $f\left(\mathbf{x}\right)$ is
congruent to the identically zero polynomial which would contradicts
our assumption that the the coefficient of the leading monomial $\prod_{0\le k<n}\left(x_{k}\right)^{-1+\left|S_{k}\right|}$
term is non-zero.$\square$

\subsection{A classical application }

Alon's Combinatorial Nullstellensatz argument is classicaly used to
prove the existence of integral coefficient Galois resolvent.\\
\textbf{Theorem} ( Integral Galois resolvent ) : For all vectors $\mathbf{r}\in\mathbb{C}$,
such that 
\begin{equation}
0\ne\prod_{0\le i<j<n}\left\langle \left(\mathbf{e}_{i}-\mathbf{e}_{j}\right),\,\mathbf{r}\right\rangle ,
\end{equation}
( where $\left\{ \mathbf{e}_{k}\right\} _{0\le k<n}$ denotes the
column vectors of the identity matrix ) there exist at least one vector
$\mathbf{s}\in\left\{ 0,\cdots,{n! \choose 2}\right\} ^{n}$ for which
the stabilizer subgroup of $S_{n}$ associated with the linear form
$\left\langle \mathbf{r},\mathbf{s}\right\rangle $ is trivial.\\
 \\
\textbf{Proof :} Consider the polynomial 
\begin{equation}
f_{\mathbf{r}}\left(\mathbf{x}\right)=\prod_{\begin{array}{c}
0\le\mbox{Lex Order}\left(\mu\right)<\mbox{Lex Order}\left(\nu\right)<n!\\
\left(\mu,\nu\right)\in S_{n}\times S_{n}
\end{array}}\left(\left\langle \mathbf{r},\,\mathbf{P}_{\mu}\mathbf{x}\right\rangle -\left\langle \mathbf{r},\,\mathbf{P}_{\nu}\mathbf{x}\right\rangle \right)
\end{equation}
 
\begin{equation}
\Rightarrow f_{\mathbf{r}}\left(\mathbf{x}\right)=\sum_{\gamma\in S_{\left(n!\right)}}\mbox{Sgn}\left(\gamma\right)\prod_{\begin{array}{c}
0\le\mbox{Lex Order}\left(\sigma\right)<n!\\
\sigma\in S_{n}
\end{array}}\left\langle \mathbf{r},\,\mathbf{P}_{\sigma}\mathbf{x}\right\rangle ^{\gamma\left(\mbox{Lex Order}\left(\sigma\right)\right)},
\end{equation}
where for an arbitrary $\sigma\in S_{n}$, the matrix $\mathbf{P}_{\sigma}$
denotes the permutation matrix 
\begin{equation}
\mathbf{P}_{\sigma}=\sum_{0\le k<n}\mathbf{e}_{k}\cdot\left(\mathbf{e}_{\sigma\left(k\right)}\right)^{T},
\end{equation}
and Lex Order($\sigma$) denotes the lexicographical order number
associated with $\sigma$. Incidentally, the leading monomial of $f_{\mathbf{r}}\left(\mathbf{x}\right)$
for some appropriately chosen monomial order is of the form 
\begin{equation}
\left(\prod_{0\le i<n}x_{i}\right)^{{n! \choose 2}},
\end{equation}
 and hence by the Combinatorial Nullstellensatz we have that 
\begin{equation}
\forall\,\left\{ A_{k}\subset\mathbb{C}\right\} _{0\le k<n},\:\mbox{ with }\left\{ \left|A_{k}\right|>{n! \choose 2}\right\} _{0\le k<n},\quad\exists\,\mathbf{a}\in A_{0}\times\cdots\times A_{n-1}\mbox{ s.t. }f_{\mathbf{r}}\left(\mathbf{a}\right)\ne0
\end{equation}
 and in particular 
\begin{equation}
\exists\,\mathbf{s}\in\left\{ 0,\cdots,{n! \choose 2}\right\} ^{n}\mbox{ s.t. }f_{\mathbf{r}}\left(\mathbf{s}\right)\ne0.
\end{equation}

\section{The Combinatorial Nullstellensatz method}

\subsection{Nullstellensatz approach to subgraph isomorphism}

Given $n\times n$ adjacency matrices $\mathbf{A}$ and $\mathbf{B}$
associated with unweighted directed graphs $G$ and $H$. We say that
$G\supseteq H$ i.e. $H$ is subisomorphic to $G$ if the following
matrix equality holds for some $n\times n$ matrix $\mathbf{P}$ 
\begin{equation}
\begin{cases}
\begin{array}{ccc}
\left(\mathbf{P}^{T}\cdot\mathbf{A}\cdot\mathbf{P}\right)\star\mathbf{B} & = & \mathbf{B}\\
\mathbf{P}^{T}\cdot\mathbf{P} & = & \mathbf{I}\\
\mathbf{P}^{\star^{2}} & = & \mathbf{P}
\end{array}\end{cases},
\end{equation}
where $\mathbf{M}\star\mathbf{N}$ denotes the entry-wise or Hadamard
product of the matrices $\mathbf{M}$, $\mathbf{N}$ and furthermore
the matrix $\mathbf{M}^{\star^{k}}$denotes the matrix resulting from
raising all non-zero entries of $\mathbf{M}$ to some integer power
$k$. We now express the matrix constraints above as constraints over
elements of the ring $\nicefrac{\mathbb{C}\left[x_{0},x_{1}\right]}{\left\{ \left(x_{j}\right)^{n}-1\right\} _{0\le j<2}}$.
We recall that the adjacency polynomials for the input graphs whose
vertices are labeled with roots of unity, are deduced from the adjacency
matrices $\mathbf{A},\mathbf{B}\in\left\{ 0,1\right\} ^{n\times n}$
as follows 
\[
A\left(x_{0},\, x_{1}\right)=n^{-2}\sum_{0\le k_{0},\, k_{1}<n}\left\langle \mathbf{w}^{\star^{-k_{0}}},\,\mathbf{w}^{\star^{-k_{1}}}\right\rangle _{\mathbf{A}}\left(x_{0}\right)^{k_{0}}\left(x_{1}\right)^{k_{1}},
\]
\[
B\left(x_{0},\, x_{1}\right)=n^{-2}\sum_{0\le k_{0},\, k_{1}<n}\left\langle \mathbf{w}^{\star^{-k_{0}}},\,\mathbf{w}^{\star^{-k_{1}}}\right\rangle _{\mathbf{B}}\left(x_{0}\right)^{k_{0}}\left(x_{1}\right)^{k_{1}},
\]
Incidentally, the subgraph isomorphism constraints are reformulated
as follows 
\[
B\left(x_{0},\, x_{1}\right)=
\]

\begin{equation}
B\left(x_{0},\, x_{1}\right)A\left(\sum_{0\le k_{0}<n}r_{k_{0}}\prod_{0\le t_{0}\ne k_{0}<n}\left(\frac{x_{0}-e^{i\frac{2\pi t_{0}}{n}}}{e^{i\frac{2\pi k_{0}}{n}}-e^{i\frac{2\pi t_{0}}{n}}}\right),\,\sum_{0\le k_{1}<n}r_{k_{1}}\prod_{0\le t_{1}\ne k_{1}<n}\left(\frac{x_{1}-e^{i\frac{2\pi t_{1}}{n}}}{e^{i\frac{2\pi k_{1}}{n}}-e^{i\frac{2\pi t_{1}}{n}}}\right)\right)
\end{equation}
\begin{equation}
\forall\;0<t<n,\quad0=\sum_{0\le j<n}\left(r_{j}\right)^{t}\mbox{ and }\: n=\sum_{0\le j<n}\left(r_{j}\right)^{n}
\end{equation}
\begin{equation}
\forall\:0\le j<2,\quad\left(x_{j}\right)^{n}=1
\end{equation}
Fortunately, the constraints may be more concisely expressed using
a polynomial parametrization of permutations of roots of unity expressed
by 
\begin{equation}
p\left(x,\mathbf{r}\right)=\sum_{0\le k<n}r_{k}\prod_{0\le s\ne k<n}\left(\frac{x-e^{i\frac{2\pi}{n}s}}{e^{i\frac{2\pi}{n}k}-e^{i\frac{2\pi}{n}s}}\right)\mod\left\{ \begin{array}{c}
\mathbf{r}^{\star^{n}}-\mathbf{w}^{\star^{0}}\\
x^{n}-1
\end{array}\right\} .
\end{equation}
The solution to the subgraph isomorphism problem is therefore completely
determined by the existence of $\gamma\in S_{n}$ and $\mathbf{g}\left(\mathbf{x},\,\mathbf{r}\right)\in\left(\mathbb{C}\left[\mathbf{x},\,\mathbf{r}\right]\right)^{n}$
such that for 
\[
\mathbf{P}_{\gamma}=\sum_{0\le k<n}\mathbf{e}_{k}\cdot\mathbf{e}_{\gamma\left(k\right)}^{T}
\]
we have 
\begin{equation}
\left\langle \left(\mathbf{r}-\mathbf{P}_{\gamma}\mathbf{w}\right),\,\mathbf{g}\left(\mathbf{x},\,\mathbf{r}\right)\right\rangle \equiv B\left(x_{0},x_{1}\right)\left[1-A\left(p\left(x_{0},\mathbf{r}\right),\, p\left(x_{1},\mathbf{r}\right)\right)\right]\mod\left\{ \begin{array}{c}
\mathbf{r}^{\star^{n}}-\mathbf{w}^{\star^{0}}\\
\mathbf{x}^{\star^{n}}-\mathbf{1}_{2\times1}
\end{array}\right\} .
\end{equation}
Since it is possible to efficiently compute the reduced polynomial
\begin{equation}
B\left(x_{0},x_{1}\right)\left[1-A\left(p\left(x_{0},\mathbf{r}\right),\, p\left(x_{1},\mathbf{r}\right)\right)\right]\mod\left\{ \begin{array}{c}
\mathbf{r}^{\star^{n}}-\mathbf{w}^{\star^{0}}\\
\mathbf{x}^{\star^{n}}-\mathbf{1}_{2\times1}
\end{array}\right\} ,
\end{equation}
it follows that it must be NP-hard to determine whether or not some
arbitrary multivariate polynomial $f\left(\mathbf{x},\,\mathbf{r}\right)\in\mathbb{C}\left[\mathbf{x},\,\mathbf{r}\right]$
admits an expansion of the form 
\begin{equation}
f\left(\mathbf{x},\,\mathbf{r}\right)=\left\langle \left(\mathbf{r}-\mathbf{P}_{\gamma}\mathbf{w}\right),\,\mathbf{g}\left(\mathbf{x},\,\mathbf{r}\right)\right\rangle 
\end{equation}
for some permutation $\gamma\in S_{n}$ and a vector $\mathbf{g}\left(\mathbf{x},\,\mathbf{r}\right)\in\left(\mathbb{C}\left[\mathbf{x},\,\mathbf{r}\right]\right)^{n}$.
\\
\\
Let Aut$f$ denote the automorphism group of $f\in\mathbb{C}\left[\mathbf{x},\,\mathbf{r}\right]$,
defined by 
\begin{equation}
\mbox{Aut }f\;:=\left\{ \sigma\in S_{n},\mbox{ s.t. }\quad f\left(\mathbf{x},\,\mathbf{r}\right)-f\left(\mathbf{x},\,\mathbf{P}_{\sigma}\mathbf{r}\right)=0\right\} ,
\end{equation}
and hence Aut$f$ denotes the stabilizer subgroup of $S_{n}$ which
fixes $f$, under permutation of the entries of the symbolic vector
$\mathbf{r}$.\\
\\
\textbf{Theorem} ( Combinatorial resolvent ) : The reduced polynomial
\begin{equation}
f\left(\mathbf{x},\mathbf{r}\right)\::=\left[1-A\left(p\left(x_{0},\mathbf{r}\right),\, p\left(x_{1},\mathbf{r}\right)\right)\right]B\left(x_{0},x_{1}\right)\mod\left\{ \begin{array}{c}
\mathbf{r}^{\star^{n}}-\mathbf{w}^{\star^{0}}\\
\mathbf{x}^{\star^{n}}-\mathbf{1}_{2\times1}
\end{array}\right\} ,
\end{equation}
admits an expansion of the form 
\begin{equation}
\left\langle \left(\mathbf{r}-\mathbf{P}_{\gamma}\mathbf{w}\right),\,\mathbf{g}\left(\mathbf{x},\,\mathbf{r}\right)\right\rangle \equiv\left[1-A\left(p\left(x_{0},\mathbf{r}\right),\, p\left(x_{1},\mathbf{r}\right)\right)\right]B\left(x_{0},x_{1}\right)\mod\left\{ \begin{array}{c}
\mathbf{r}^{\star^{n}}-\mathbf{w}^{\star^{0}}\\
\mathbf{x}^{\star^{n}}-\mathbf{1}_{2\times1}
\end{array}\right\} .
\end{equation}
for some permutation $\gamma\in S_{n}$ and a vector $\mathbf{g}\left(\mathbf{x},\,\mathbf{r}\right)\in\left(\mathbb{C}\left[\mathbf{x},\,\mathbf{r}\right]\right)^{n}$,
if and only if
\begin{equation}
0\equiv\prod_{\sigma\in\nicefrac{S_{n}}{\mbox{Aut}f}}f\left(\mathbf{x},\,\mathbf{P}_{\sigma}\mathbf{r}\right)\mod\left\{ \begin{array}{c}
\mathbf{r}-\mathbf{w}\\
\mathbf{x}^{\star^{n}}-\mathbf{1}_{2\times1}
\end{array}\right\} 
\end{equation}
 \\
\textbf{Proof} : The proof of the theorem is an immediate consequence
of Euclidean division. We have 
\begin{equation}
\forall\,\sigma^{-1}\in S_{n},\; f\left(\mathbf{x},\,\mathbf{r}\right)=\kappa_{\sigma^{-1}}\left(\mathbf{x}\right)+\left\langle \left(\mathbf{r}-\mathbf{P}_{\sigma^{-1}}\mathbf{w}\right),\,\mathbf{g}_{\sigma^{-1}}\left(\mathbf{x},\,\mathbf{r}\right)\right\rangle 
\end{equation}
\begin{equation}
\Rightarrow f\left(\mathbf{x},\,\mathbf{P}_{\sigma^{-1}}\mathbf{r}\right)=\kappa_{\sigma^{-1}}\left(\mathbf{x}\right)+\left\langle \mathbf{P}_{\sigma^{-1}}\left(\mathbf{r}-\mathbf{w}\right),\,\mathbf{g}_{\sigma^{-1}}\left(\mathbf{x},\,\mathbf{P}_{\sigma^{-1}}\mathbf{r}\right)\right\rangle 
\end{equation}
\begin{equation}
\Rightarrow f\left(\mathbf{x},\,\mathbf{P}_{\sigma^{-1}}\mathbf{r}\right)=\kappa_{\sigma^{-1}}\left(\mathbf{x}\right)+\left\langle \left(\mathbf{r}-\mathbf{w}\right),\:\mathbf{P}_{\sigma}\mathbf{g}_{\sigma^{-1}}\left(\mathbf{x},\,\mathbf{P}_{\sigma^{-1}}\mathbf{r}\right)\right\rangle 
\end{equation}
and hence 
\[
\prod_{\sigma\in S_{n}}f\left(\mathbf{x},\,\mathbf{P}_{\sigma^{-1}}\mathbf{r}\right)\equiv\prod_{\sigma\in S_{n}}\kappa_{\sigma^{-1}}\left(\mathbf{x}\right)\mod\left\{ \begin{array}{c}
\mathbf{r}-\mathbf{w}\\
\mathbf{x}^{\star^{n}}-\mathbf{1}_{2\times1}
\end{array}\right\} 
\]
furthermore we note that 
\begin{equation}
\prod_{\sigma\in S_{n}}f\left(\mathbf{x},\,\mathbf{P}_{\sigma^{-1}}\mathbf{r}\right)\equiv\left(\prod_{\sigma\in\nicefrac{S_{n}}{\mbox{Aut}f}}f\left(\mathbf{x},\,\mathbf{P}_{\sigma}\mathbf{r}\right)\right)^{\left|\mbox{Aut}f\right|}
\end{equation}
from which it immediately follows that 
\begin{equation}
\prod_{\sigma\in\nicefrac{S_{n}}{\mbox{Aut}f}}f\left(\mathbf{x},\,\mathbf{P}_{\sigma}\mathbf{r}\right)\equiv0\mod\left\{ \begin{array}{c}
\mathbf{r}-\mathbf{w}\\
\mathbf{x}^{\star^{n}}-\mathbf{1}_{2\times1}
\end{array}\right\} \Leftrightarrow\exists\:\sigma\in S_{n}\:\mbox{ s.t. }\kappa_{\sigma}\left(\mathbf{x}\right)\equiv0.\:\square
\end{equation}

\subsection{Subgraph isomorphism dual Nullstellensatz construction}

In order to mimic the Alon and Tarsi polynomial constructions discussed
in \cite{AT} for determining the existence of solutions to a subgraph
isomorphism instance, one would seek instead a polynomial construction
which will be identically zero if $H$ is not sub-isomorphic to $G$,
and the polynomial construction would admit a non vanishing term otherwise.
For the purpose of the construction let us consider the adjacency
polynomials associated with the graphs 
\[
A\left(x_{0},\, x_{1}\right)=n^{-2}\sum_{0\le k_{0},\, k_{1}<n}\left\langle \mathbf{w}^{\star^{-k_{0}}},\,\mathbf{w}^{\star^{-k_{1}}}\right\rangle _{\mathbf{A}}\:\left(x_{0}\right)^{k_{0}}\left(x_{1}\right)^{k_{1}},
\]
\[
B\left(x_{0},\, x_{1}\right)=n^{-2}\sum_{0\le k_{0},\, k_{1}<n}\left\langle \mathbf{w}^{\star^{-k_{0}}},\,\mathbf{w}^{\star^{-k_{1}}}\right\rangle _{\mathbf{B}}\:\left(x_{0}\right)^{k_{0}}\left(x_{1}\right)^{k_{1}},
\]
deduced as usual from the adjacency matrices $\mathbf{A}$, $\mathbf{B}\in\left\{ 0,1\right\} ^{n\times n}$.
Furthermore consider the set $\mathfrak{G}_{\mathbf{B}}$ which denote
the set of adjacency matrices of non isomorphic graphs which do not
contain $H$ as a subgraph. Formally we write 
\begin{equation}
\mathfrak{G}_{\mathbf{B}}\,:=\left\{ \mathbf{C}\in\left\{ 0,\,1\right\} ^{n\times n},\:\mbox{ s.t. }\begin{array}{c}
\forall\sigma\in S_{n}\,\left(\mathbf{1}_{n\times n}-\mathbf{P}_{\sigma}^{T}\mathbf{C}\mathbf{P}_{\sigma}\right)\star\mathbf{B}\ne\mathbf{0}_{n\times n}\\
\forall\sigma\in S_{n}\mbox{ and}\,\left(\mathbf{C}_{0},\mathbf{C}_{1}\right)\in\left(\mathfrak{G}_{B}\backslash\left\{ \mathbf{C}_{1}\right\} \right)\times\left(\mathfrak{G}_{B},\backslash\left\{ \mathbf{C}_{0}\right\} \right),\:\left(\mathbf{P}_{\sigma}^{T}\mathbf{C}_{0}\mathbf{P}_{\sigma}\right)\ne\mathbf{C}_{1}
\end{array}\right\} 
\end{equation}
finally let $\mathfrak{P}_{\mathbf{B}}$ denote the corresponding
set of adjacency polynomials defined as 
\begin{equation}
\mathfrak{P}_{\mathbf{B}}:=\left\{ n^{-2}\sum_{0\le k_{0},\, k_{1}<n}\left\langle \mathbf{w}^{\star^{-k_{0}}},\,\mathbf{w}^{\star^{-k_{1}}}\right\rangle _{\mathbf{C}}\:\left(x_{0}\right)^{k_{0}}\left(x_{1}\right)^{k_{1}},\;\mbox{ s.t. }\mathbf{C}\in\mathfrak{G}_{\mathbf{B}}\right\} .
\end{equation}
The corresponding Alon and Tarsi polynomial construction for determining
the existence of solution to subgraph Isomorphism instances is expressed
by 
\[
f_{B}\left(A,\,\mathbf{x},\,\mathbf{r}\right)=\prod_{0\le i<j<n}\left\langle \left(\mathbf{e}_{i}-\mathbf{e}_{j}\right),\,\mathbf{r}\right\rangle \times
\]
\begin{equation}
\prod_{C\left(x_{0},x_{1}\right)\in\mathfrak{P}_{\mathbf{B}}}\left[\left[1-A\left(x_{0},x_{1}\right)\right]C\left(p\left(x_{0},\mathbf{r}\right),p\left(x_{1},\mathbf{r}\right)\right)\right]\mod\left\{ \begin{array}{c}
\mathbf{r}^{\star^{n}}-\mathbf{w}^{\star^{0}}\\
\mathbf{x}^{\star^{n}}-\mathbf{1}_{2\times1}
\end{array}\right\} .
\end{equation}
It follows by construction that if $f_{B}\left(A,\,\mathbf{x},\,\mathbf{r}\right)$
is identically zero then this fact constitutes a certificate of non-existence
of solutions to the subisomorphism instance. Incidentally, we think
of the construction as dual to the initial polynomial encoding of
subgraph isomorphism problem described in the previous section. Similarly,
the dual polynomial construction for graph Isomorphism instances is
expressed by 
\[
g_{B}\left(A,\,\mathbf{x},\,\mathbf{r}\right)=\prod_{0\le i<j<n}\left\langle \left(\mathbf{e}_{i}-\mathbf{e}_{j}\right),\,\mathbf{r}\right\rangle \times
\]
\begin{equation}
\prod_{C\left(x_{0},x_{1}\right)\in\mathfrak{P}_{\mathbf{B}}}\left[\left[1-A\left(x_{0},x_{1}\right)\right]C\left(p\left(x_{0},\mathbf{r}\right),p\left(x_{1},\mathbf{r}\right)\right)+\left[1-C\left(p\left(x_{0},\mathbf{r}\right),p\left(x_{1},\mathbf{r}\right)\right)\right]A\left(x_{0},x_{1}\right)\right]\mod\left\{ \begin{array}{c}
\mathbf{r}^{\star^{n}}-\mathbf{w}^{\star^{0}}\\
\mathbf{x}^{\star^{n}}-\mathbf{1}_{2\times1}
\end{array}\right\} .
\end{equation}

\section*{Acknowledgments}

This material is based upon work supported by the National Science
Foundation under agreements Princeton University Prime Award No. CCF-0832797
and Sub-contract No. 00001583. The author would like to thank the
IAS for providing excellent working conditions. The author is also
grateful to Vladimir Retakh, Ahmed Elgammal, Avi Wigderson, Noga Alon
for insightful comments while preparing this manuscript.


\begin{thebibliography}{References}
\bibitem[A]{A}N. Alon, Combinatorial Nullstellensatz, Comb. Prob.
Comput. 8 , 7-29, (1999).

\bibitem[AT]{AT} N. Alon and M. Tarsi Colorings and orientations
of graphs. Combinatorica 12 125\textendash{}134, (1992) 

\bibitem[AKS]{AKS} M. Agrawal, N. Kayal and N. Saxena, PRIMES is
in P, Ann. Math. v. 160, 781\textendash{}793, (2004)

\bibitem[B]{B} Bruno Buchberger. An algorithmic criterion for the
solvability of a system of algebraic equations. Aequationes Mathematicae
4, 374\textendash{}383, (1970).

\bibitem[CLO]{CLO} D. A. Cox J. B. L.D. O'Shea. Ideals, Varieties,
and Algorithms Third Edition, Springer, (2007).

\bibitem[H]{H}P. C. Heinig, Proof of the combinatorial nullstellensatz
over integral domains in the spirit of Kouba, Electron. J. Combin.
17 (2010).

\bibitem[K]{K}E. Kranakis: Invited Talk: Symmetry and Computability
in Anonymous Networks. SIROCCO 1-16, (1996).

\bibitem[Ko]{Ko}J. Kollár, ``Sharp Effective Nullstellensatz'',
J. of AMS. 1 (4): 963\textendash{}975 (1988).

\bibitem[Kou]{Kou}O. Kouba, A Duality Based Proof of the Combinatorial
Nullstellensatz, Elect. J. of Combinatorics, V. 16, (2009).

\bibitem[La]{La} M. Laso\'{n}, A generalization of Combinatorial
Nullstellensatz, Electron. J. Combin. 17 (2010).

\bibitem[LHMO]{LHMO}J. A. Loera, C. J. Hillar, P. N. Malkin, M. Omar,
Recognizing Graph Theoretic Properties with Polynomial Ideals, Elect.
J. of Combinatorics, V. 17, (2010).

\bibitem[LMM]{LMM}J. A. De Loera, J. Lee, P. N. Malkin, and S. Margulies.
Hilbert's nullstellensatz and an algorithm for proving combinatorial
infeasibility. In Proceedings of the twenty-first international symposium
on Symbolic and algebraic computation (ISSAC 2008). ACM, New York,
NY, USA, 197-206. (2008).

\bibitem[LMO]{LMO}J. A. Loera, J. Lee, S. Margulies, and S. Onn.
Expressing combinatorial problems by systems of polynomial equations
and hilbert's nullstellensatz. Comb. Probab. Comput. 18, 4 , 551-582
( 2009).

\bibitem[Lo]{Lo}L. Lovász. Stable sets and polynomials. Discrete
Mathematics, 124:137\textendash{}153, (1994).

\bibitem[M]{M}S. Margulies, Computer Algebra, Combinatorics, and
Complexity: Hilbert's Nullstellensatz and NP-Complete Problems. Ph.D.
Dissertation. University of California at Davis, Davis, CA, USA. AAI3336295.
(2008)

\bibitem[Mi]{Mi}M. Michalek, A short proof of Combinatorial Nullstellensatz,
Amer. Math. Monthly 117 , 821\textendash{}823, (2010).

\bibitem[RW]{RW}R. Ramamurthi, D. B. West, Hypergraph extension of
the Alon\textendash{}Tarsi list coloring theorem, Combinatorica 25
, 355-366 (2005).

\bibitem[S]{S}U. Schauz, Algebraically solvable problems: describing
polynomials as equivalent to explicit solutions. Electron. J. Combin.
15 , no. 1, Research Paper 10 (2008).

\bibitem[SSK]{SSK}S. Ben-Israel, Eli Ben-Sasson, David R. Karger:
Breaking local symmetries can dramatically reduce the length of propositional
refutations. Electronic Colloquium on Computational Complexity (ECCC)
17: 68 (2010)\end{thebibliography}
\end{document}